\documentclass{amsart}
\usepackage{palatino,hhline, bbm}
\usepackage{amsthm, amsmath}
\usepackage{amssymb} 
\usepackage{amscd}
\usepackage{mathrsfs}
\renewcommand{\mathcal}{\mathscr}
\usepackage[hmargin=3cm, vmargin=3cm]{geometry}
\usepackage[breaklinks,bookmarksopen,bookmarksnumbered]{hyperref}
\usepackage[all]{xy}

\theoremstyle{plain}

\newtheorem*{lem}{Lemma}
\newtheorem{prop}{Proposition}
\newtheorem*{cor}{Corollary}

\theoremstyle{remark}
\newtheorem*{rem}{Remark}

\newcommand\pr{\noindent\textit{Proof} : }

\newcommand\rond{\kern 1pt{\scriptstyle\circ}\kern 1pt}

\newcommand\Pic{\operatorname{Pic}}

\newcommand\Z{\mathbb{Z}}

\newcommand\C{\mathbb{C}}
\renewcommand\P{\mathbb{P}}

\newcommand\G{\mathbb{G}}
\renewcommand\O{\mathcal{O}}

\newcommand\iso{\vbox{\hbox to .8cm{\hfill{$\scriptstyle\sim$}\hfill}
\nointerlineskip\hbox to .8cm{{\hfill$\longrightarrow $\hfill}} }}
\newcommand\bir{\vbox{\hbox to .8cm{\hfill{$\scriptstyle\sim$}\hfill}
\nointerlineskip\hbox to .8cm{{\hfill$\dasharrow $\hfill}} }}
\newcommand\abs[1]{\lvert {#1}\rvert}

\begin{document}
\title{An ampleness criterion for rank 2 vector bundles on surfaces}
\author[Arnaud Beauville]{Arnaud Beauville}
\address{Universit\'e C\^ote d'Azur\\
CNRS -- Laboratoire J.-A. Dieudonn\'e\\
Parc Valrose\\
F-06108 Nice cedex 2, France}
\email{arnaud.beauville@unice.fr}
 
\maketitle 
\section{Introduction}
We observe in this note that  the proof of the Bogomolov stable restriction theorem \cite{B} can be adapted to give
an  ampleness criterion for  globally generated rank 2 vector bundles on certain surfaces. This applies to the Lazarsfeld-Mukai bundles, to  congruences of lines in  $\P^3$, and possibly to the construction of surfaces with ample cotangent bundle. 

\medskip	
\section{Main result}	
Throughout the note, $S$ will be a smooth projective surface over $\C$. We denote by $N^1(S)$ the group of divisors on $S$ modulo numerical equivalence; this is a free, finitely generated abelian group,  quotient of $\operatorname{NS}(S)=H^2(S,\Z)_{\mathrm{alg}} $ by its torsion subgroup. 

\begin{prop}\label{fib}
Let $E$ be a globally generated rank 2 vector bundle on $S$, with $h^0(E)\geq 4$. Assume that $N^1(S)=\Z\cdot c_1(E)$. Then either $E$ is ample, or $E=\O_S\oplus  \det (E)$.
\end{prop}
We will need the following lemma:

\begin{lem}
Let $S$ be a smooth projective surface, and let $E$ be a globally generated rank 2 vector bundle on $S$, with $h^0(E)\geq 4$ and $H^1(S,\det(E)^{-1})=0$. Then $c_1^2(E)>c_2(E)$. 
\end{lem}
\pr Let $V$ be a general $4$-dimensional subspace of $H^0(S,E)$. Then $V$ generates $E$ globally, giving rise to an exact sequence 
\begin{equation}\label{se}
0\rightarrow N\rightarrow V\otimes _{\C}\O_S\rightarrow E\rightarrow 0\, .
\end{equation}
 Since $N^*$ is globally generated, the zero locus of a general section $s$ of $N^*$ is finite, of length  $c_2(N^*)=\allowbreak c_1^2(E)-c_2(E)$. Thus this number is $\geq 0$; if it is zero, $s$ does not vanish, so we have  an exact sequence
\[0\rightarrow \O_S \xrightarrow{\ s\ } N^* \rightarrow \det(E)\rightarrow 0\,.\]
Since $H^1(S,\det(E)^{-1})=0$, this sequence splits, so that $N\cong \O_S\oplus \det(E)^{-1}$. Thus the exact sequence (\ref{se}) reduces to
\[0\rightarrow \det(E)^{-1} \rightarrow \O_S^3\rightarrow E\rightarrow 0\,;\]
but using again $H^1(S,\det(E)^{-1})=0$ this implies  $h^0(E)\leq 3$, contradicting the hypothesis.\qed

\medskip	
\noindent
\emph{Proof of the Proposition} : We denote  by $c_1$ and $c_2$ the Chern classes of $E$ in $H^*(S,\Z)$, and by $\Delta _E:=4c_2-c_1^2$ its discriminant. Assume that $E$ is not ample. By Gieseker's lemma \cite[Proposition 6.1.7]{L}, there exists an irreducible curve $C$ in $S$ and a surjective homomorphism $u:E\twoheadrightarrow \O_C$.   The kernel $F$ of $u$ is a vector bundle, with total Chern class
$c(F)=c(E)c(\O_C)^{-1}=\allowbreak (1+c_1+c_2)(1-[C])$, hence 
\[c_1(F)=c_1-[C]\ ,\ c_2(F)=c_2-c_1\cdot [C] \ ,\ \mbox{and }\ \Delta _F=\Delta _E-2c_1\cdot [C]-[C]^2\,.\]
  
The curve $C$ is numerically equivalent to $rc_1$ for some integer  $r\geq 1$. Therefore  \[\Delta _F=4c_2 -(r+1)^2c_1^2\leq 4(c_2-c_1^2)\,.\]Because of our hypotheses   $\det(E)$ is ample, so $H^1(S,\det(E)^{-1}) =0$ and we can apply the Lemma, which gives $\Delta _F<0$. By Bogomolov's theorem (see \cite[Th\'eor\`eme 6.1]{Ra}),   we have an exact sequence
\[0\rightarrow L \rightarrow F \rightarrow \mathcal{I}_ZM\rightarrow 0\]
where $Z$ is a finite subscheme of $S$, $L$ and $M$ are line bundles on $S$, with $c_1(L)=ac_1$,   $c_1(M)=bc_1$ for some integers $a,b$ such that $a\geq b$.

From that exact sequence we get $c_1(F)=(a+b)c_1$, hence $a+b=1-r$, and $c_2(F)=\deg(Z)+ abc_1^2$, hence $\Delta _F=\allowbreak 4\deg(Z)-(a-b)^2c_1^2$. Comparing with the previous expression for $\Delta _F$ and using the Lemma  again we find 
\[(a-b)^2c_1^2\geq -\Delta _F=(r+1)^2c_1^2-4c_2> (r^2+2r-3)c_1^2\geq (r^2-1)c_1^2\ ,\]
hence  $a-b\geq r$, and $a\geq 1$. 

We have $H^0(E\otimes L^{-1})=H^0(E^*\otimes \det(E)\otimes L^{-1})\neq 0$. Since $E$ is globally generated, the natural map $E^*\rightarrow H^0(E)^*\otimes _{\C}\O_S$ is injective, hence
 $H^0(\det(E)\otimes L^{-1})\neq 0$. Since $c_1(L)=ac_1$ with $a\geq 1$, the only possibility is  $L\cong \det(E)$,   and therefore $H^0(E^*)\neq 0$. Using again that $E$ is globally generated, we obtain 
 $E=\O_S\oplus \det(E)$.\qed

\begin{rem}
The condition $h^0(E)\geq 4$ is necessary: if $E$ is ample and globally generated, the rational map $\P(E)\rightarrow \P(H^0(E))$ associated to the linear system $\abs{\O_{\P(E)}(1)}$ is a finite morphism, hence $\dim \P(H^0(E))\geq 3$. On the other hand, the condition $N^1(S)=\Z\cdot c_1$ is quite restrictive, but it is not clear how it could be weakened. For instance, we will exhibit in Example 1 of \S4 a globally generated rank 2 vector bundle $E$ on $\P^2$ with $h^0(E)\geq 4$, $\det E=\O_{\P^2}(2)$, which is not ample. 
\end{rem}

\section{Application 1: Lazarsfeld-Mukai bundles}

Let $C$ be an irreducible  curve in $S$,  $L$  a line bundle on $C$, and   $V$  a 2-dimensional subspace of $H^0(L)$ which generates $L$. The \emph{Lazarsfeld-Mukai bundle} $E_{C,V}$ is defined by the exact sequence
\[0\rightarrow E^*_{C,V}\rightarrow V\otimes_{\C} \O_S \rightarrow L\rightarrow 0\,.\]
Let $N_C:=\O_S(C)_{|C}$ be the normal of $C$ in $S$. By duality we get an exact sequence
\[0\rightarrow V^*\otimes_{\C} \O_S \rightarrow E_{C,V}\rightarrow N_C\otimes L^{-1}\rightarrow 0\,.\]
\begin{prop}
Assume   $H^1(S,\O_S)=0$,   $N^1(S)=\Z\cdot [C]$, and that the line bundle $N_C\otimes L^{-1}$ on $C$ is globally generated and nontrivial. Then $E_{C,V}$ is globally generated and ample.
\end{prop}

\pr We put $E:=E_{C,V}$. Since  $H^1(S,\O_S)=0$, we have a commutative diagram of exact sequences
\[\xymatrix{0\ar[r] & V^*\otimes _{\C}\O_S\ar[r]\ar@{=}[d] & H^0(S,E) \otimes_{\C}\O_S\ar[r]\ar[d] & H^0(C,N_C\otimes L^{-1})\otimes _{\C}\O_S\ar[r]\ar@{->>}[d] & 0\\
0\ar[r] & V^*\otimes _{\C}\O_S\ar[r]  & E \ar[r] & N_C\otimes L^{-1}\ar[r]& 0\,.
}\]This implies that $E$ is globally generated, with $h^0(E)=2+h^0(N_C\otimes L^{-1})\geq 4$. From the bottom exact sequence we get $c_1(E)=[C]$ and
 $c_2(E)= \deg(L) >0$. The conclusion follows from Proposition~\ref{fib}.\qed

\section{Application 2: congruences of lines}

\smallskip	

Let $\G$ be the Grassmannian of lines in $\P^3$, which we view as a smooth quadric in $\P^5$; let $S\subset \G$ be a smooth surface. This defines a 2-dimensional family of lines in $\P^3$, classically called a \emph{congruence}. A point $p\in\P^3$ through which pass infinitely many lines of the congruence is called a \emph{fundamental point} (or, more classically, a singular point) of the congruence. 

\begin{prop}
Assume that $S$ has degree $>1$ and that $N^1(S)$ is generated by the restriction of $\O_{\G}(1)$. Then $S$ has no fundamental point. 
\end{prop}
\pr Let $E$ be the restriction to $S$ of the universal quotient bundle $Q$ on $\G$. The projective bundle $\P(E)$ on $S$ parametrizes pairs $(\ell,p)$ in $S\times \P^3$ with $p\in\ell$, and the second projection  $q:\P_S(E)\rightarrow \P^3$ satisfies $q^*\O_{\P^3}(1)=\O_{\P(E)}(1)$.
 Thus $q$ is finite (that is, $S$ has no fundamental point) if and only if $E$ is ample.

We have $h^0(Q)=4$, and a nonzero section of $Q$ vanishes along a linear plane; therefore $h^0(E)\geq 4$, and we can apply Proposition \ref{fib}. If $E=\O_S\oplus \O_S(1)$, we have $c_2(E)=0$, that is, $c_2(Q)\cdot [S]=0$; this can only happen if $S$ is a linear plane, which we have excluded. Therefore $E$ is ample. \qed

\begin{cor}
Let $d,e$ be two integers with  $d,e>1$, or $d=1$ and $e\geq 3$; let $S\subset \G$ be the complete intersection of two general hypersurfaces of degree $d$ and $e$. Then $S$ has no fundamental point.
\end{cor}
Indeed $\Pic(S)$ is generated by $\O_S(1)$ \cite[Th\'eor\`eme 1.2]{D}.

\medskip	
\noindent\emph{Examples}$.-$ 1) Perhaps the simplest example of a nontrivial congruence is the surface $S$ of  lines bisecant to a twisted cubic $T\subset\P^3$; it is isomorphic to $\operatorname{Sym}^2T\cong \P^2 $, embedded in $\G\subset\P^5$ by the Veronese map. In that case $N^1=\Z\cdot [\O_{\P^2}(1)]$ but $\det E=\O_{\P^2}(2)$, and indeed the  fundamental locus of $S$ is $T$, so $E$ is not ample.

\smallskip	

2) Let $A$ be an abelian surface such that $\operatorname{NS}(A)=\Z\cdot [L] $, where $L$ is a line bundle with $L^2=10$. The linear system $\abs{L}$ embeds $A$ into $\P^4$ \cite{R}, giving the famous Horrocks-Mumford abelian surface. The projection $\pi :\G\rightarrow \P^4$ from a general point of $\P^5$ is a double covering, and the surface   $S:=\pi ^{-1}(A)\subset\G$ is smooth. The line bundle $\pi ^*L$ is not divisible in $N^1(S)$: since 
$(\pi ^*L)^2=20$, this could happen only if $\pi ^*L$ is divisible by 2; but $\pi ^*L=K_S$, so this would imply that $K_S^2$ is divisible by 8, a contradiction. It then follows from \cite{Bu} that $N^1(S)$ is generated by $\pi ^*L=\O_S(1)$, so Proposition 2 applies and $S$ has no fundamental point.

\bigskip	
\section{Application 3 (virtual): surfaces with ample cotangent bundle}

\smallskip	
The original motivation of this work was to obtain new examples of surfaces with ample cotangent bundle -- these surfaces have very interesting properties, but there are few concrete examples known. Applying Proposition \ref{fib} to $\Omega ^1_S$ we get the following result;  unfortunately we do not know any example of a surface satisfying the hypotheses 
(help welcome!).
\begin{prop}
Assume that $\Omega ^1_S$ is globally generated (for instance that $S$ is a subvariety of an abelian variety), $q(S)\geq 4$, and $N^1(S)=\Z\cdot [K_S]$. Then $\Omega ^1_S$ is ample.
\end{prop}
\pr The hypotheses imply that $K_S$ is ample, hence $c_2(S)>0$; therefore $\Omega ^1_S$ is not isomorphic to $\O_S\oplus K_S$. The conclusion follows from
 Proposition \ref{fib}. \qed


\bigskip	
\begin{thebibliography}{X-}

\bibitem[B]{B} F. Bogomolov\,: \textsl{Holomorphic tensors and vector bundles on projective varieties}. Math. of the USSR, Izvestija \textbf{13}  (1979), 499-555.

\bibitem[Bu]{Bu} A. Buium\,: \textsl{Sur le nombre de Picard des rev\^etements doubles des surfaces alg\'ebriques}. 
C. R. Acad. Sci. Paris S\'er. I Math. \textbf{296}  (1983), no. 8, 361-364.

\bibitem[D]{D} P. Deligne\,: \textsl{Le th\'eor\`eme de Noether}. Groupes de monodromie en g\'eom\'etrie alg\'ebrique II, Expos\'e XIX, 328-340. Lecture Notes in Math. \textbf{340} , Springer, Berlin-Heidelberg-New York, 1973.


\bibitem[L]{L} R. Lazarsfeld\,: \textsl{Positivity in algebraic geometry}, II. Ergeb.  Math. (3) \textbf{49}.  Springer-Verlag, Berlin, 2004.

\bibitem[R]{R} S. Ramanan\,: \textsl{Ample divisors on abelian surfaces}. 
Proc. London Math. Soc. (3) \textbf{51}  (1985), no. 2, 231-245. 


\bibitem[Ra]{Ra} M. Raynaud\,: \textsl{Fibr\'es vectoriels instables -- applications aux surfaces (d'apr\`es Bogomolov)}.  Algebraic surfaces, 293-314; Lecture Notes in Math. \textbf{868} , Springer, Berlin-New York, 1981.


\end{thebibliography}
\end{document}